\documentclass{amsart}

\usepackage{amssymb}
\usepackage{amsthm}

\title[Isomorphism problem for Abelian $p$-Groups]{The isomorphism problem for computable Abelian $p$-groups of 
bounded length}


\author{Wesley Calvert}

\address{Department of Mathematics\\
255 Hurley Hall\\
University of Notre Dame\\
Notre Dame, Indiana, 46556}

\email{wcalvert@nd.edu}

\thanks{The author was partially supported by NSF Grants DMS 9970452 and DMS 0139626. The author wishes to thank J.\
F.\ Knight for many helpful comments on this paper.}

\newtheorem{thm}{Theorem}[section]
\newtheorem{prop}[thm]{Proposition}
\newtheorem{lem}[thm]{Lemma}

\newtheorem{claim}[thm]{Claim}

\theoremstyle{definition}
\newtheorem{dfn}[thm]{Definition}

\newtheorem{Q}[thm]{Question}

\begin{document}

\begin{abstract} Theories of classification distinguish classes with some good structure theorem from those for which
none is possible. Some classes (dense linear orders, for instance) are non-classifiable in general, but are
classifiable when we consider only countable members. This paper explores such a notion for classes of computable
structures by working out a sequence of examples. 

We follow recent work by Goncharov and Knight in using the degree of the isomorphism problem for a class to
distinguish classifiable classes from non-classifiable.  In this paper, we calculate the degree of the isomorphism 
problem for Abelian $p$-groups of bounded Ulm length.  The result is a sequence of classes whose isomorphism problems 
are cofinal in the hyperarithmetical hierarchy.  In the process, new back-and-forth relations on such groups are 
calculated. \end{abstract}

\maketitle

\section{Introduction}

In an earlier paper \cite{ecl1}, we began to consider a notion of ``classification" for classes of computable
structures.  For some classes, there is a ``classification," or ``structure theorem" of some kind.  For instance, the
classification of algebraically closed fields states that a single cardinal (the transcendence degree) completely
determines the structure up to isomorphism.  For other classes (graphs, for example, or arbitrary groups) such a
result would be surprising, and when we introduce the necessary rigor we can prove that there is none to be found.  
They simply have more diversity than any structure theorem could describe.

We assume all structures have for a universe some computable subset of $\omega$ and identify a structure with its atomic diagram.  Thus, for instance, a structure is
computable if and only if its atomic diagram is computable, as a set of G\"{o}del numbers of sentences.  Alternatively, we could use the quantifier-free diagram instead of
the atomic diagram.  Similarly, a structure is associated with the index of a Turing machine which enumerates its atomic diagram (assuming its universe is computable).  In
this paper, I will write $\mathcal{A}_a$ for the computable structure with atomic diagram $W_a$ and will always assume that a class $K$ of structures has only computable
members.  The following definition was recently proposed by Goncharov and Knight \cite{gk}.

\begin{dfn}The isomorphism problem, denoted $E(K)$, is the set \[\{ (a,b) | \mathcal{A}_a , \mathcal{A}_b \in K\mbox{, and }\mathcal{A}_a \simeq 
\mathcal{A}_b\}\]\end{dfn}

If the set of indices for computable members of $K$, denoted $I(K)$, is hyperarithmetical, then $E(K)$ is $\Sigma^1_1$.  Intuitively, in the worst case, where $E(K)$ is
properly $\Sigma^1_1$, the easiest way to say that two members of $K$ are isomorphic is to say, ``There exists a function which is an isomorphism between them."  Often
there are easier ways to check isomorphism, such as counting basis elements of vector spaces.  Such a ``shortcut" is a classification.  There is also a natural ``floor" to 
the complexity of $E(K)$, since to say that $a \in I(K)$ requires saying that $a$ is an index for some structure, which is already $\Pi^0_2$.  

This notion is closely related to work in descriptive set theory, originating in the work of Friedman and Stanley
\cite{fs}.  In that context, the set of countable models of a theory is viewed as a topological space, and we would
calculate the topological complexity of the isomorphism relation as a subset of the Cartesian product of two copies of
the space (for a more complete description of the topological situation, see \cite{hjbook}).  Many of the proofs that
a class has maximal complexity, like the Friedman -- Stanley proof of the Borel completeness of fields \cite{fs},
require only minor modification.

Several classes are well-known to have maximally complicated isomorphism problems.  The following theorem summarizes several classical results.  Proofs may be found in 
articles by Rabin and Scott \cite{rabsco}, Goncharov and Knight \cite{gk}, Morozov \cite{moro}, and Nies \cite{nies}.
\begin{thm} \label{cla} If $K$ is the set of computable members of any of the following classes, then $E(K)$ is 
$\Sigma^1_1$ complete:
\begin{enumerate}
\item Undirected graphs
\item Linear orders
\item Trees
\item Boolean algebras
\item Abelian $p$-groups
\end{enumerate}
\end{thm}

The following additions to the list follow easily from recent work by Hirschfeldt, Khoussainov, Shore, and Slinko \cite{hkss}.
\begin{thm} [Hirschfeldt -- Khoussainov -- Shore -- Slinko] \label{hksst}If $K$ is the set of computable members of 
any 
of the following classes, then $E(K)$ is $\Sigma^1_1$ complete:
\begin{enumerate}
\item Rings
\item Distributive lattices
\item Nilpotent groups
\item Semigroups
\end{enumerate}
\end{thm}

In an earlier paper \cite{ecl1}, the following were added:
\begin{thm} \label{ecf} \rule{0in}{0in} \newline
\vspace{-.2in}
\begin{enumerate}
\item If $K$ is the set of computable members of any of the following classes, then $E(K)$ is $\Sigma^1_1$ complete:
  \begin{enumerate}
   \item Fields of any fixed characteristic
   \item Real Closed Fields
  \end{enumerate}
\item If $K$ is the set of computable members of any of the following classes, then $E(K)$ is $\Pi^0_3$ complete:
  \begin{enumerate}
   \item Vector spaces over a fixed computable field
   \item Algebraically closed fields of fixed characteristic
   \item Archimedean real closed fields
  \end{enumerate}
\end{enumerate}
\end{thm}

In this paper, the complexity of the isomorphism problem will be calculated for other classes.  Two major goals which are partially achieved here are the
answers to the following questions:
\begin{Q} What are the possible complexities of the isomorphism problem for classes of structures? \end{Q} 
\begin{Q} Do classes with high complexity acquire it all at once?\end{Q}
Considering Abelian $p$-groups of bounded Ulm length will give us a sequence of isomorphism problems whose degrees are
cofinal in the hyperarithmetical degrees.  In some sense, this also shows a smooth transition from very low complexity
(say, $\Pi^0_3$ complete) to the ``non-classifiable" (that is, properly $\Sigma^1_1$).

\section{Notation and Terminology for Abelian $p$-groups} 

Let $p$ be an arbitrary prime number.  Abelian $p$-groups are Abelian groups in which each element has some power of
$p$ for its order.  We will consider only countable Abelian $p$-groups.  These groups are of particular interest
because of their classification up to isomorphism by Ulm.  For a classical discussion of this theorem and a more
detailed discussion of this class of groups, consult Kaplansky's book \cite{kaplansky}.  Generally, notation here will
be similar to Kaplansky's.

It is often helpful to follow L.\ Rogers \cite{lrogers} in representing these groups by trees.  Consider a tree $T$.  The Abelian $p$-group $G(T)$ is the group generated by
the nodes in $T$ (among which the root is $0$), subject to the relations stating that the group is Abelian and that $px$ is the predecessor of $x$ in the tree.  Reduced
Abelian $p$-groups, from this perspective, are represented by trees with no infinite paths.

The idea of Ulm's theorem is that it generalizes the notion that to determine a finitely generated torsion Abelian
group it is only necessary to determine how many cyclic components of each order are included in a direct sum
decomposition. Let $G$ be an Abelian $p$-group.  We will produce an ordinal sequence (usually transfinite) of
cardinals $u_{\beta}(G)$ (each at most countable), which is constant after some ordinal (called the ``length" of $G$).  
If $H$ is also an Abelian $p$-group and for all $\beta$ we have $u_{\beta} (G) = u_{\beta} (H)$, then $H \simeq G$
(this is still subject to another condition we have yet to define).

First set $G_0 = G$.  Now we inductively define $G_{\beta + 1} = pG_{\beta} = \{ px | x \in G_{\beta}\}$, where $px$
denotes the sum of $x$ with itself $p$ times. We also define, for limit $\beta$, the subgroup $G_{\beta} =
\bigcap\limits_{\gamma < \beta} G_{\gamma}$.  Further, let $P(G)$ denote the subgroup of elements $x$ for which $px =
0$, and let $P_{\beta} (G) = P \cap G_{\beta}$.  Now the quotient $P_{\beta} (G)/ P_{\beta + 1}(G)$ is a
$\mathbb{Z}_p$ vector space, and we call its dimension $u_{\beta} (G)$.  Where no confusion is likely, we will omit
the argument $G$ and simply write $P_{\beta}$, and so forth.

For any Abelian $p$-group $G$, there will be some least ordinal $\lambda (G)$ such that $G_{\lambda (G)} = G_{\lambda (G) + 1}$.  This is called the \textit{length} of $G$.  
If $G_{\lambda (G)} = \{0\}$, then we say that $G$ is \textit{reduced}.  Equivalently, $G$ is reduced if and only if 
it has no divisible subgroup. The \textit{height} of an
element $x$ is the unique $\beta$ such that $x \in G_{\beta}$, but $x \notin G_{\beta + 1}$.  It is conventional to write $h(0) = \infty$, where $\infty$ is greater than
any ordinal.  Similarly, if our group contains a divisible element $x$, we write $h(x) = \infty$.  In the course 
of this paper, we will only consider reduced groups.  When $G$ is a direct sum of cyclic groups, 
$u_n (G)$ is exactly equal to the number of direct summands of order $p^{n+1}$.  We can now state Ulm's theorem,
but we will not prove it here.

\begin{thm} [Ulm] Let $G$ and $H$ be reduced countable Abelian $p$-groups.  Then $G \simeq H$ if and only if for every
countable ordinal $\beta$ we have $u_{\beta}(G)~=~u_{\beta}(H)$. \end{thm}

It is interesting to note that this theorem is not ``recursively true."  Lin showed that if two computable groups
satisfying the hypotheses of this theorem have identical Ulm invariants, they may not be computably isomorphic
\cite{clin}.  However, it is known that (depending heavily on the particular statement of the theorem), Ulm's theorem
is equivalent to the formal system $\mbox{ATR}_0$ \cite{frsimsm, simpson}.  Related work from a constructivist
perspective may be found in a paper by Richman \cite{richman}.

A calculation of the complexity of the isomorphism problem for special classes of computable reduced Abelian $p$-groups is essentially a computation of the complexity of 
checking the equality of Ulm invariants.  Given some computable ordinal $\alpha$, we will consider the class of reduced Abelian $p$-groups of length at most $\alpha$.

\section{Bounds on Isomorphism Problems}
When we begin to consider special classes of Abelian $p$-groups from the perspective described in section 1, it
quickly becomes apparent that all examples in Theorems \ref{cla}, \ref{hksst}, and \ref{ecf} were especially nice 
ones.  In all of these cases,
$I(K)$ was $\Pi^0_2$ and $E(K)$ was something worse.  It is easy to see that $I(K) \leq_T E(K)$, since
$$I(K) = \{a | (a, a) \in E(K)\}$$

For instance, if $K$ is the class of reduced Abelian $p$-groups of length at most $\omega$, $I(K)$ is $\Pi^0_3$ complete.  Then to show that $E(K)$ is $\Pi^0_3$ complete, 
it is enough to show that $E(K)$ is $\Pi^0_3$, and this is not difficult (the reader interested in the details of this may wish to glance ahead to Proposition 
\ref{komega}).

However, this doesn't tell us whether $E(K)$ has high complexity ``on its own,"  or just by virtue of it being hard to tell whether we have something in $K$.  In a talk in
Almaty in the summer of 2002, J.\ Knight proposed the following definition to clear up the distinction: \begin{dfn} Suppose $A \subseteq B$.  Let $\Gamma$ be some
complexity class (e.g.\ $\Pi^0_3$), and $K$ a class of computable structures.  Then $A$ is $\Gamma$ \emph{within} $B$ 
if and only if there is some $R \in \Gamma$ such that $A = R
\cap B$ \end{dfn}

In the example above, saying that $E(K)$ is $\Pi^0_3$ within $I(K) \times I(K)$ means that there is a $\Pi^0_3$
relation $R(a,b)$ such that if $a$ and $b$ are indices for computable reduced Abelian $p$-groups, then $R(a,b)$
defines the relation ``$\mathcal{A}_a$ has the same Ulm invariants as $\mathcal{A}_b$."  In general, it is possible
that $A$ is not $\Gamma$ but that $A$ is $\Gamma$ \textit{within} $B$.  Consider for instance the case of a theory
which is $\aleph_0$-categorical.  If $K$ is the class of models of such a theory, then $E(K)$ is not computable, but
$E(K)$ is computable within $I(K) \times I(K)$.

We can also define a reducibility ``within $B$", which will, in turn, give us a notion of completeness.

\begin{dfn} Let $A, B$, and $\Gamma$ be as in the previous definition.
\begin{enumerate} 
\item $S \leq_m A$ \emph{within} $B$ if there is a computable $f: \omega \to B$ such that for all $n$, $n \in S \iff 
f(n) \in A$.
\item $A$ is $\Gamma$ \emph{complete within} $B$ if $A$ is $\Gamma$ within $B$ and for any $S \in \Gamma$ we have $S 
\leq_m A$ within $B$.
\end{enumerate} \end{dfn}

Essentially, this definition says that $A$ is $\Gamma$ complete within $B$ if it is $\Gamma$ within $B$ and there is a function witnessing that it is $\Gamma$ complete
which only calls for questions about things in $B$.  In fact, the questions are only about members of a c.e.\ subset of $B$.  We will usually write ``within $K$" for
``within $I(K) \times I(K)$."  All results stated in section 1 remain true when we add ``within $K$"  to their statements, and the original proofs still work.  In fact,
this is intuitively the ``right" way to say that the structure of a class is complicated: we say that if we look at some members, it is difficult to tell whether they are
isomorphic.  It would be unconvincing to argue that the structure of a class is complicated simply because it is difficult to tell whether things are in the class or not.

For any computable ordinal $\alpha$, it is somewhat straightforward to write a computable infinitary sentence stating that $G$ is a reduced Abelian $p$-group of length at
most $\alpha$ and that $G$ and $H$ have the same Ulm invariants up to $\alpha$.  In particular, Barker \cite{barker} verified the following.

\begin{lem} Let $G$ be a computable Abelian $p$-group.
\begin{enumerate}
\item $G_{\omega \cdot \alpha}$ is $\Pi^0_{2 \alpha}$.
\item $G_{\omega \cdot \alpha + m}$ is $\Sigma^0_{2 \alpha + 1}$.
\item $P_{\omega \cdot \alpha}$ is $\Pi^0_{2 \alpha}$.
\item $P_{\omega \cdot \alpha + m}$ is $\Sigma^0_{2 \alpha + 1}$.
\end{enumerate}
\end{lem}

\begin{proof} It is easy to see that 3 and 4 follow from 1 and 2 respectively.  Toward 1 and 2,
note the following: 
\begin{eqnarray*}
x \in G_m & \iff & \exists y (p^m y = x)\\
x \in G_{\omega} & \iff & \bigwedge\limits_{m \in \omega} \hspace{-0.15in}\bigwedge \exists y (p^m y = x)\\
x \in G_{\omega \cdot \alpha + m} & \iff & \exists y [p^m y = x \wedge G_{\omega \cdot \alpha}(y)]\\
x \in G_{\omega \cdot \alpha + \omega} & \iff & \bigwedge\limits_{m \in \omega} \hspace{-0.15in}\bigwedge \exists y [p^m y = x \wedge G_{\omega \cdot \alpha}(y)]\\
x \in G_{\omega \cdot \alpha} & \iff & \bigwedge\limits_{\gamma < \alpha} \hspace{-0.15in}\bigwedge G_{\omega \cdot \gamma} (x) \mbox{ for limit $\alpha$}\\
\end{eqnarray*}
\end{proof}

Work by Lin \cite{linjsl}, when viewed from our perspective, shows that for any $m \in \omega$, there is a group $G$ in which $G_m$ is $\Sigma^0_1$ complete.  Given this
lemma, we can place bounds on the complexity of $I(K)$ and $E(K)$.

\begin{lem} If $K_{\alpha}$ is the class of reduced Abelian $p$-groups of length at most $\alpha$, then $I(K_{\omega \cdot \beta + m})$ is $\Pi^0_{2 
\beta + 1}$.\end{lem}

\begin{proof}
The class $K_{\omega \cdot m + \beta}$ may be characterized by the axioms of an Abelian $p$-group (which are 
$\Pi^0_2$), together with the condition $$\forall x [x \in G_{\omega \cdot 
\beta +m} \rightarrow x = 0]$$  Since the previous lemma guarantees that this condition is $\Pi^0_{2 \beta +1}$, we 
know that $I(K_{\omega \cdot \beta +m})$ is also $\Pi^0_{2 \beta +1}$.
\end{proof}

\begin{lem} \label{ebnd} If $K_\alpha$ is as in the previous lemma, we use $\hat{\alpha}$ to denote $\mbox{$\sup\limits_{\omega \cdot \gamma < \alpha} (2 \gamma + 3)$}$.  
Then $E(K_{\alpha})$ is $\Pi^0_{\hat{\alpha}}$ within $K$.\end{lem}

\begin{proof}
Note that the relation ``there are at least $n$ elements of height $\beta$ which are $\mathbb{Z}_p$-independent over $G_{\beta + 1}$" is
defined in the following way.  To say that $x_i, \dots, x_n$ are $\mathbb{Z}_p$-independent over $G_{\beta+1}$, we write the computable $\Pi^0_{2 \beta + 1}$ formula 
\[D_{n, \beta}(x_1, \dots, x_n) = \bigwedge\limits_{b_1, \dots b_n \in \mathbb{Z}_p} (\sum\limits_{i = 1}^n b_i x_i \notin G_{\beta+1})\]
Now to write ``there are at least $n$ independent elements of height $\beta$ and order $p$," we use the sentence
\[B_{n, \beta} = \exists x_1, \dots, x_n [(\bigwedge\limits_{i = 1}^n G_\beta (x_i)) \wedge (\bigwedge\limits_{i=1}^n px_i = 0)) \wedge D_{n, \beta} (\overline{x})]\]
which is a computable $\Sigma^0_{2 \beta + 2}$ sentence.  Now we can define isomorphism by
\[\bigwedge\limits_{\rule{0.05in}{0in} n\in\omega \atop  \rule{0.05in}{0in}\beta<\alpha}\hspace{-0.18in}\bigwedge \mathcal{A}_a \models B_{n,\beta} \Leftrightarrow
\mathcal{A}_b \models B_{n,\beta}\] 
We write each $\beta < \alpha$ as $\beta = \omega \cdot \gamma + m$, where $m \in \omega$.  If $\hat{\alpha}$ is as defined in the statement of the lemma, then this can be expressed by a computable $\Pi^0_{\hat{\alpha}}$ sentence.
\end{proof}

\section{Completeness for Length $\omega \cdot m$}

\begin{prop} \label{komega} If $K_\omega$ is the class of computable Abelian $p$-groups of length at most $\omega$, then $E(K_{\omega})$ is $\Pi^0_3$ complete within
$K_\omega$.\end{prop}

\begin{proof} We first observe that the set is $\Pi^0_3$ within $K$, by applying the previous lemma.  Now let $S = \forall e \exists \tilde{y} \forall z
\overline{R}(n,e,y,z)$ be an arbitrary $\Pi^0_3$ set.  We can represent $S$ as the set defined by \[\forall e \exists^{<\infty} y \hspace{0.05in}R(n,e,y)\] where
$\exists^{<\infty}$ is read ``there exist at most finitely many."  Consider the Abelian $p$-group $G^{\omega}$ with Ulm sequence
\[u_\alpha = \left\{\begin{array}{ll}
    \omega & \mbox{if $\alpha < \omega$}\\
    0 & \mbox{otherwise}\\ \end{array} \right. \]

We will build a uniformly computable sequence $H^n$ of reduced Abelian $p$~-~groups of height at most $\omega$ such
that $H^n \simeq G^{\omega}$ if and only if $n \in S$.  Let $G^{\omega,\infty}$ denote the direct sum of countably
many copies of the smallest divisible Abelian $p$-group $\mathbb{Z}(p^\infty)$, and note that $G^{\omega, \infty}$ has
a computable copy, as a direct sum of copies of a subgroup of $\mathbb{Q} / \mathbb{Z}$.  We will denote the element
where $x$ occurs in the $i$th place with zeros elsewhere by $(x)_i$.  For instance, $G^{\omega, \infty}$ set-wise is 
the collection of all sequences of proper fractions whose denominators are powers of $p$, and the element
$(\frac{1}{p})_2$ denotes the element $(0, \frac{1}{p}, 0, 0, \dots)$.

List the atomic sentences by $\phi_e$, the pairs of elements in $G^{\omega, \infty}$ by $\xi_e$, and set $D_{-1} = C_{-1} = Y_{e, -1} = X_{e, -1} = \tilde{X}_{e,-1} = T_{e, 
-1} = \emptyset$.  We will build groups to meet the
following requirements: \smallskip

\begin{tabular}{rl} 
$P_e$ : & There are infinitely many independent elements $x \in H^n$ of order $p$ and\\
        & height exactly $e$ if and only if there are at most finitely many $y$\\
        & such that $R(n,e,y)$.\\
$Q_e$ : & If $\xi_e = (a,b)$ and $a, b \in H^n$, then $a+b \in H^n$.\\
$Z_e$ : & If all parameters occurring in $\phi_e$ are in $H^n$, then exactly one of\\
        & $\phi_e \in D$ or $\lnot \phi_e \in D$.\\
\end{tabular}
\smallskip

Roughly speaking, $D_s$ will be the diagram of $H^n$, and $C_s$ will be its domain.  For each $e$, the set $Y_{e,s}$ will keep track of the $y$ already seen, $X_{e,s}$ the 
$x$ created of height at least $e$, and $\tilde{X}_{e,s}$ the $x$ which are given greater height, as in $P_e$.  The set $T_{e,s}$ will keep track of the heights greater 
than $e$ already used to put elements from $X_e$ in $\tilde{X}_e$, so that we do not accidentally make infinitely many elements of height $e+1$.

We say that $P_e$ requires attention at stage $s$ if there is some $y<s$ such that $y \notin Y_{e, s-1}$ and $R(n,e,y)$ and there is also some $x \in X_{e,s-1} \setminus
\tilde{X}_{e,s-1}$, or if for all $y <s$ we have either $y \in Y_{e,s-1}$ or $\lnot R(n,e,y)$.  We say that $Q_e$ requires attention at stage $s$ if $\xi_e = (a,b)$ and $a,
b \in C_{s-1}$ but $a+b \notin C_{s-1}$.  We say that $Z_e$ requires attention at stage $s$ if all parameters that occur in $\phi_e$ are in $C_{s-1}$ and $D_{s-1}$ does not
include either $\phi_e$ or $\lnot \phi_e$.

At stage $s$, to satisfy $P_e$, we will act by first looking for some $y<s$ such that $y \notin Y_{e, s-1}$ and 
$R(n,e,y)$.  If none is found, the action will be to enumerate a
new independent $x$ of height at least $e$.  To do this, find the first $k$ such that $(\frac{1}{p})_{k}$ does not occur in $C_{s-1}$ or in any element of $D_{s-1}$.  Let
\[C_{s} = C_{s-1} \cup \{(\frac{1}{p^j})_{k} | j = 1, \dots, (e-1) \}\] and set $X_{e,s} = X_{e, s-1} \cup \{(\frac{1}{p})_k\}$, $\tilde{X}_{e, s} = \tilde{X}_{e,s-1}$,
$T_{e,s} = T_{e,s-1}$, and $Y_{e,s} = Y_{e,s-1}$.  If such a $y$ is found, on the other hand, the action will be to 
give all existing element of $X_{e,s-1}$ height greater than
$e$.  To do this, collect \[K = \{k | (\frac{1}{p})_k \in X_{e, s-1} \setminus \tilde{X}_{e,s-1}\}\] 
and the 
least positive $r \notin T_{e,s-1}$.  Note that $K$ is finite.  Set \[C_s =
C_{s-1} \cup \bigcup\limits_{k \in K} \{(\frac{1}{p^j})_{k} | j = (e, \dots, e+r+1)\}\] and set $T_{e, s} = T_{e, s-1} 
\cup \{r\}$, $\tilde{X}_{e, s} = \tilde{X}_{e, s-1} \cup \{(\frac{1}{p})_k | k \in K
\}$, $X_{e, s} = X_{e, s-1}$, and $Y_{e,s} = Y_{e,s-1} \cup \{y\}$.

To satisfy $Q_e$ at stage $s$ we will look to see whether the elements of $\xi_e = (a, b)$ are in $C_{s-1}$.  
If they are both there, set $C_{s} = C_{s-1} \cup 
\{a+b\}$.  Otherwise, set $C_s = C_{s-1}$.

To satisfy $Z_e$, we will act at stage $s$ by first looking for the parameters in $\phi_e$ in $C_{s-1}$.  If all of 
them are there and $G^{\omega, \infty} \models \phi_e$, then
set $D_s = D_{s-1} \cup \{\phi_e\}$.  If all of them are there and $G^{\omega, \infty} \models \lnot \phi_e$, then set $D_s = D_{s-1} \cup \{\lnot \phi_e\}$.  If some of
the parameters are not in $C_{s-1}$, we set $D_s = D_{s-1}$.

Now if $n \in S$, for each $e$ we have $Q_e$ to guarantee that $u_e(H^n)$ will be infinite, so $H^n \simeq G^\omega$.  If $n \notin S$, there is some $e$ such that $Q_e$ 
guarantees that $u_e(H^n)$ is finite, so $H^n \not\simeq G^\omega$. \end{proof}

Since this result is perfectly uniform, we can use it for induction.  What we actually have established is the following:
\begin{prop} If $S$ is a set which is $\Pi^0_3$ relative to $X$, then there is a uniformly $X$-computable sequence of reduced Abelian $p$-groups $(H^n)_{n \in \omega}$,
each of length at most $\omega$, such that $H^n \simeq G^{\omega}$ if and only if $n \in S$.\end{prop}
There is a result of Khisamiev \cite{khis}, which allows us to transfer these \\ $X$-computable groups down to the
computable level.
\begin{prop}[Khisamiev] If $G$ is a $X^{\prime \prime}$-computable reduced Abelian \\ $p$-group, then there is an $X$-computable reduced Abelian
$p$-group $H$ such that $H_{\omega} \simeq G$ and $u_n(H) = \omega$ for all $n \in \omega$.  Moreover, from an index for $G$, we can effectively compute an index for 
$H$.\end{prop} 
These two results together can be used to establish

\begin{prop} \label{finite} If $K_{\omega \cdot m}$ is the class of computable reduced Abelian $p$-groups of length at most $\omega \cdot m$, then $E(K_{\omega \cdot m})$
is $\Pi^0_{2m+1}$ complete within $K$.\end{prop}

\begin{proof} Let $S$ be an arbitrary $\Pi^0_{2m+1}$ set.  Since $S$ is $\Pi^0_3$ in $\emptyset^{(2m-1)}$, we have a uniformly $\emptyset^{(2m+1)}$-computable sequence of
reduced Abelian $p$-groups $(H^n)_{n \in \omega}$, each of length at most $\omega$, such that $H^n \simeq G^{\omega}$ if and only if $n \in S$.  Now we can step each $H^n$
down to a lower level using Khisamiev's result, so that we have a uniformly $\emptyset^{(2n-3)} = \emptyset^{(2(n-1)-1)}$-computable sequence $(H^{2,n})_{n \in \omega}$ of
reduced Abelian $p$-groups, each of height $\omega \cdot 2$ which again have the property that $H^{2,n}$ has a constantly infinite Ulm sequence if and only if $n \in S$.  
By induction, we define $(H^{i,n})_{n \in \omega}$, and when we get to $(H^{m,n})_{n \in \omega}$, it will be a uniformly computable sequence of groups of length at most
$\omega \cdot m$ such that $H^{m,n}$ has constantly infinite Ulm sequence if and only if $n \in S$. \end{proof}

\section{Completeness for Higher Bounds on Length}
Giving completeness results for higher levels requires more elaborate machinery.  We will prove a more general result using an $\alpha$-system, in the sense of Ash.  These
systems are explained in detail, along with several other variants, in the book of Ash and Knight \cite{akbook}.  The ``metatheorem" for $\alpha$-systems was proved in a
paper by Ash \cite{ashlab}.

Roughly speaking, an $\alpha$-system describes all possible priority constructions of a given kind, and the metatheorem states that given an ``instruction function" which
is $\Delta^0_\alpha$, the system will produce a c.e.\ set (in our case, the diagram of a group) which incorporates the information given in the instruction
function.  More formally, we make the following definition:

\begin{dfn} [Ash] Let $\alpha$ be a computable ordinal.  An $\alpha$-system is a structure $$(L, U, P, \hat{\ell}, E, (\leq_\beta)_{\beta < \alpha})$$ where $L$ and $U$ are
c.e.\ sets, $E$ is a partial computable function on $L$ (it will eventually enumerate the diagram of the structure we are building), $P$ is a c.e.\ alternating tree on $L$
and $U$ (that is, a set of strings with letters alternating between $L$ and $U$) in which all members start with $\hat{\ell} \in L$, and $\leq_\beta$ are uniformly c.e.\
binary relations on $L$, where the following properties are satisfied:
\begin{enumerate} 
\item $\leq_\beta$ is reflexive and transitive for all $\beta < \alpha$ 
\item $a \leq_{\gamma} b \Rightarrow a \leq_{\beta} b$ for all $\beta < \gamma < \alpha$
\item If $a \leq_0 b$, then $E(a) \subseteq E(b)$ 
\item If $\sigma u \in P$, where $\sigma$ ends in $\ell^0$, and $$\ell^0 \leq_{\beta_0} \ell^1 \leq_{\beta_1} \dots 
\leq_{\beta_{k-1}} \ell^k$$ where $\beta_0 > \beta_1 >
\dots > \beta_k$, then there exists some $\ell^*$ such that $\sigma u \ell^* \in P$ and for all $i \leq k$, we have $\ell^i \leq_{\beta_i} \ell^*$.
\end{enumerate}\end{dfn}

If we have such a system, we say that an \textit{instruction function} for $P$ is a function $q$ from the set of sequences in $P$ of odd length (i.e.\ those with a last
term in $L$) to $U$, so that for any $\sigma$ in the domain of $q$, $\sigma q(\sigma) \in P$.  The following theorem, due to Ash \cite{ashlab}, guarantees that if we have
such a function, there is a string which represents ``carrying out" the instructions while enumerating a c.e.\ set.  
We call an infinite string $\pi = \hat{\ell}u_1 \ell_1 u_2 \ell_2 \dots$ a ``run" of
$(P, q)$ if it is a path through $P$ with the property that for any initial segment $\sigma u$ we have $u = q(\sigma)$.  The metatheorem also guarantees that
there is a run with the property that $\bigcup\limits_{i \in \omega} E(\ell_i)$ is computably enumerable.

\begin{prop} [Ash Metatheorem] If we have an $\alpha$-system $$(L, U, P, \hat{\ell}, E, (\leq_{\beta})_{\beta < \alpha})$$ and if $q$ is a $\Delta^0_\alpha$ instruction
function for $P$, then there is a run $\pi: \omega \to (L \cup U)$ of $(P,q)$ such that $\bigcup\limits_{i \in \omega} E(\pi(2i))$ is c.e.  Further, from computable indices
for the components of the system and a $\Delta^0_\alpha$ index for $q$, we can effectively determine a c.e.\ index for $\bigcup\limits_{i \in \omega} E(\pi(2i))$.\end{prop}

What this means is that if we can set up an appropriate system, then given some highly undecidable requirements, we can build a computable group to satisfy them.  The
difficulty (aside from digesting the metatheorem itself) mainly consists of defining the right system.  Afterwards, it is no trouble to write out the high-level
requirements we want to meet.  Using such a system, we will prove the following generalization of Proposition \ref{finite}.

\begin{thm} Let $\alpha$ be a computable limit ordinal, and let $\hat{\alpha} = \sup\limits_{\omega \cdot \gamma < \alpha} (2 \gamma + 3)$, as in Proposition \ref{ebnd}.  
If $K_{\alpha}$ is the class of reduced Abelian $p$-groups of length at most $\alpha$ then $E(K_{\alpha})$ is $\Pi^0_{\hat{\alpha}}$ complete within $K_\alpha$.\end{thm}

\begin{proof} Let $(\alpha_i)_{i \in \omega \setminus \{0\}}$ be a sequence cofinal in $\alpha$ (for instance, if $\alpha = \omega \cdot \omega$, then $\alpha_i = \omega
\cdot i$ would do, or if $\alpha = \omega \cdot (\beta +1)$, we could use $\alpha_i = \omega \cdot \beta + i$).  Consider the family of groups $(\hat{G}^i)_{i \in \omega}$,
each of length $\alpha$ where $\hat{G}_0$ has uniformly infinite Ulm sequence and
$$u_\beta (\hat{G}^i) =
\left\{\begin{array}{ll} \omega &
 \mbox{if $\beta < \alpha_i$ or if $\beta$ is even}\\ 0 & \mbox{otherwise}\\
 \end{array} \right.$$ 
Since the Ulm sequences of these groups are uniformly computable, there is a uniformly computable sequence $(G^i)_{i \in \omega}$ such that $G^i \simeq \hat{G}^i$ for all
$i$, and such that in each of these groups, for any $\beta$, the predicate ``$x$ has height $\beta$" is computable.  The proof of this, which is due to Oates, is a
modification of an argument of L.\ Rogers \cite{lrogers}, and may be found in Barker's paper \cite{barker}.

For any set $S \in \Pi^0_{\hat{\alpha}}$, we will construct a sequence of groups $(H^n)_{n \in \omega}$ such that if $n \in S$ then $H^n \simeq G^0$, and otherwise, $H^n
\simeq G^i$ for some $i \neq 0$.  To do this, we will define an $\hat{\alpha}$-system.  Let $L$ be the set of pairs $(j, p)$, where $j \in \omega$ and $p$ is a finite
injective partial function from $\omega$ to $G^j$.  Let $U$ be the set $\{0,1\}$.  By $E(j, p)$, we will mean the first $|dom(p)|$ atomic or negation atomic sentences
with parameters from the image of $p$ which are true in $G^j$.  Let
$\hat{\ell} = (0, \emptyset)$, and $P$ be the set of strings of the form $\hat{\ell} u_1 \ell_1 u_2 \ell_2 \dots$ which satisfy the following properties:
\begin{enumerate}
\item $u_i \in U$ and $\ell_i \in L$
\item If $u_i = 1$ then $u_{i+1} = 1$
\item If $\ell_i = (j_i , p_i)$, then both the domain and range of $p_i$ contain at least the first $i$ members of $\omega$
\item If $\ell_i = (j, p)$ and $u_i = 1$, then $j \neq 0$.  Otherwise, $j=0$.  Further, if $u_{i-1} = 1$ and $\ell_{i-1} = (j_{i-1}, q)$, then $j = j_{i-1}$.
\end{enumerate}

For the $\leq_{\beta}$ we will modify the standard back-and-forth relations on Abelian $p$-groups.  In general, the standard back-and-forth relations on a class $K$ are
characterized as relations on pairs $(\mathcal{A}, \overline{a})$ where $\mathcal{A} \in K$ and $\overline{a}$ is a finite tuple of $\mathcal{A}$.  

\begin{dfn} If $\overline{a} \subseteq \mathcal{A}$ and $\overline{b} \subseteq \mathcal{B}$ are finite tuples of 
equal length, then we define the \emph{standard back-and-forth relations} $\leq_\beta$ as follows:
\begin{enumerate}
\item $(\mathcal{A}, \overline{a}) \leq_1 (\mathcal{B}, \overline{b})$ if and only if for all finitary $\Sigma^0_1$ formulas true of $\overline{b}$ in $\mathcal{B}$ are 
true of $\overline{a}$ in $\mathcal{A}$.
\item $(\mathcal{A}, \overline{a}) \leq_\beta (\mathcal{B}, \overline{b})$ if and only if for any finite $\overline{d} \subset \mathcal{B}$ and any $\gamma$ with $1 \leq 
\gamma < \beta$ there is some $\overline{c} \subset \mathcal{A}$ of equal length such that $(\mathcal{B}, \overline{b}, \overline{d}) \leq_\gamma (\mathcal{A}, 
\overline{a}, \overline{c})$.
\end{enumerate}
\end{dfn}

This definition extends naturally to tuples of different length as follows: we say that 
$(\mathcal{A},\overline{a}) \leq_\beta (\mathcal{B}, \overline{b})$ if and only if
$\overline{a}$ is no longer than $\overline{b}$ and that for the initial segment $\overline{b}^\prime \subset \overline{b}$ of length equal to that of $\overline{a}$, we
have $(\mathcal{A}, \overline{a}) \leq_\beta (\mathcal{B}, \overline{b}^\prime)$.  Barker \cite{barker} gave a useful characterization of these relations in the case of
Abelian $p$-groups $\mathcal{A}$ and $\mathcal{B}$, where $\mathcal{A} = \mathcal{B}$.

\begin{prop}[Barker] If $\leq_{\beta}$ are the standard back-and-forth relations on reduced Abelian $p$-groups, and if
$\overline{a}$ and $\overline{b}$ are finite subsets of equal length in an Abelian $p$-group with the height of
elements given by $h$ respectively and with equal cardinality, with a function $f$ mapping elements of $\bar{b}$ to 
corresponding elements of $\bar{a}$, then the following hold:
\begin{enumerate}
\item $\overline{a} \leq_{2 \cdot \delta} \overline{b}$ if and only if the two generate isomorphic subgroups and for every $b \in \overline{b}$ and $a = f(b)$ we have 
$$h(a) =
h(b) < \omega \cdot \delta \mbox{ or } h(b), h(a) \geq \omega \cdot \delta$$
\item $\overline{a} \leq_{2 \cdot \delta + 1} \overline{b}$ if and only if the two generate isomorphic subgroups and for every $b \in \overline{b}$ and $a = f(b)$ we have
\begin{enumerate}
\item In the case that $P_{\omega \cdot \delta + k}$ is infinite for every $k \in \omega$, $$h(a) = h(b) < \omega \cdot \delta$$ or $$h (b) \geq \omega \cdot 
\delta \mbox{ and } h (a) \geq \min \{h (b), \omega \cdot \delta + \omega\}$$
\item In the case that $P_{\omega \cdot \delta + k}$ is infinite and $P_{\omega \cdot \delta + k+1}$ is finite, 
$$h(a) = h(b) < \omega \cdot \delta$$ or 
$$\omega \cdot \delta \leq h(b) \leq h (a) \leq \omega \cdot \delta + k$$ or 
$$h(a) = h(b) > \omega \cdot \delta + k$$
\item In the case that $P_{\omega \cdot \delta}$ is finite, $$h(x) = h(x)$$
\end{enumerate}
\end{enumerate}
\end{prop}

Since in all groups with which we are concerned, $P_{\omega \cdot \delta + k}$ will be infinite for all $\delta < \alpha$, we will have no need for the more complicated 
cases.  Also, it is helpful to deal with groups which satisfy the stronger condition that they 
have infinite Ulm invariants at each limit level.

\begin{dfn} Let $\mathcal{A}, \mathcal{B}$ be countable reduced Abelian $p$-groups of length at most
 $\alpha$ such that for any limit ordinal $\nu < \alpha$ we have $u_{\nu}(\mathcal{A}) 
= u_{\nu}(\mathcal{B}) = \omega$.  Let the height of an element in
its respective group be given by $h$.  Let $\overline{a}, \overline{b}$ be finite sequences of 
equal length from $\mathcal{A}$ and $\mathcal{B}$, respectively.  Then
define $(\leq_{\delta})_{\delta < \omega_1}$ by the following:
\begin{enumerate}
\item $(\mathcal{A}, \overline{a}) \leq_{2 \cdot \delta} (\mathcal{B}, \overline{b})$ if and only if 
\begin{enumerate}
\item The function matching elements of $\overline{a}$ to corresponding elements of $\overline{b}$ extends to an 
isomorphism $f: <\overline{b}> \to <\overline{a}>$, 
\item for every $b \in \overline{b}$ and $a = f(b)$ we have $$h(a) = h(b) < \omega \cdot \delta \mbox{ or } h(b), h(a) \geq \omega \cdot \delta$$ and
\item for all $\beta < \omega \cdot \delta$ we have $u_\beta (\mathcal{A}) = u_{\beta} (\mathcal{B})$.
\end{enumerate}
\item $(\mathcal{A}, \overline{a}) \leq_{2 \cdot \delta + 1} (\mathcal{B}, \overline{b})$ if and only if 
\begin{enumerate}
\item The function matching respective elements in $\overline{a}$ and $\overline{b}$ extends to an isomorphism $f: <\overline{b}> \to <\overline{a}>$, 
\item for every $b \in \overline{b}$ and $a = f(b)$ we have $$h(a) = h(b) < \omega \cdot \delta$$ or $$h (b) \geq \omega \cdot \delta \mbox{ and } h (a) \geq \min \{h (b),
\omega \cdot \delta + \omega\}$$
\item for all $\beta < \omega \cdot \delta$ we have $u_\beta (\mathcal{A}) = u_{\beta} (\mathcal{B})$.
\item for all $\beta \in [\omega \cdot \delta, \omega \cdot \delta + \omega)$ we have $u_\beta (\mathcal{A}) \geq u_{\beta} (\mathcal{B})$.
\end{enumerate}
\end{enumerate}
\end{dfn}

In order to verify that we have an $\hat{\alpha}$-system, the following lemma will be important.

\begin{lem} \label{baf} Suppose $(\mathcal{A}, \overline{a}) \leq_{\beta} (\mathcal{B}, \overline{b})$.  Then for any $\eta < \beta$ and for any finite sequence
$\overline{d} \subseteq \mathcal{B}$ there exists a sequence $\overline{c} \subseteq \mathcal{A}$ of equal length such that $(\mathcal{B}, \overline{b}, \overline{d})
\leq_{\eta} (\mathcal{A}, \overline{a}, \overline{c})$.\end{lem}

\begin{proof} Suppose that the conditions stated for $\leq_{2 \cdot \delta}$ hold.  Now suppose $\delta~=~\gamma~+~1$.  
It suffices to show that for all finite sequences
$\overline{d} \subseteq \mathcal{B}$ there exists a sequence $\overline{c} \subseteq \mathcal{A}$ of equal length such that $(\mathcal{B}, \overline{b}, \overline{d})
\leq_{2 \cdot \delta + 1} (\mathcal{A}, \overline{a}, \overline{c})$.  We will extend $f$ to $\overline{d}$ one element at a time.  Let $d \in \overline{d}$, and suppose
that $d \notin <\overline{b}>$ (since if it were in that subgroup, we could simply map it to the corresponding element of $<\overline{c}>$.  Further suppose, without loss
of generality, that $pd \in <\overline{b}>$ and that $h(d) \geq h(d + s)$ for any $s \in <\overline{b}>$.  This last condition is often stated ``$d$ is proper with respect
to $<\overline{b}>$."  These assumptions are reasonable, since if we need to extend $f$ to an element farther afield, we can go one element at a time and work down to it.  
From this point, we essentially follow Kaplansky's proof of Ulm's theorem \cite{kaplansky} to find the appropriate match for $d$.  Use $z$ to denote $f(pd)$.  It now
suffices to find some $c$ of height $h(d)$ which is proper with respect to $<\overline{a}>$ and such that $pc = z$.

First suppose that $h(z) = h(d) + 1$.  Now both $z$ and $pd$ must be nonzero.  For $c$ we may choose any element of $(\mathcal{A})_{h(d)}$ with $pc = z$.  The height of
$z$ tells us that there must exist such an element.  We first check that $h(c) \leq h(d)$, which is easy, since if $h(c) > h(d)$, we would have \[h(z) = h(pc) \geq h(c)+1
\gneq h(d)  + 1\]  Finally, it is necessary to show that $c$ is proper with respect to $<\overline{a}>$.  Suppose that $c \in <\overline{a}>$.  Then $c = f(y)$ for some $y
\in <\overline{b}>$.  Then $pd = py$ and $d-y \notin <\overline{b}>$ to avoid $d \in <\overline{b}>$.  Further, $h(d-y) = h(d)$, since $h(y) = h(d)$ and $d$ is proper with
respect to $<\overline{b}>$.  However, \[h(p(x-y)) = h(0) = \infty \gneq h(d) + 1\] contradicting the maximality of 
$h(px)$.  Thus $c \notin <\overline{a}>$.  Now suppose we have
$h~(c~+~t)~\geq~h(d)~+~1$ for some $r \in <\overline{a}>$ with $r = f(s)$.  Since $c+r \neq 0$ (to avoid the case 
that
$c =-r \in <\overline{a}>$), we know that $h(p(w+r)) \geq h(d) + 2$,
so that $h(p(d+s)) \geq h(d) + 2$.  Since $h(r) \geq h(d)$, we also have $h(s) \geq h(d)$, so $h(d+s) = h(d)$, contradicting the maximality of $h(pd)$.

Suppose that $h(z) > h(d) + 1$.  Now there is some $v \in (\mathcal{B})_{h(d) + 1}$ such that $pd = pv$.  Then the element $d-v$ is in $P_{h(d)}(\mathcal{B})$, has 
height $h(d)$, and is thus proper with respect to $<\overline{b}>$.  I make the following claim.

\begin{claim} [Lemma 13 of \cite{kaplansky}] Let the function $$r: (<\overline{b}>_{h(d)} \cap p^{-1} (\mathcal{B})_{h(d+2)}) \to P_{h(d)}(\mathcal{B})$$
 be defined as follows: For any $x \in (<\overline{b}>_{h(d)} \cap p^{-1} (\mathcal{B})_{h(d)+2})$ there exists some $y \in (\mathcal{B})_{h(d)+1}$ such that $py=px$.  
Define $Y$ by $Y: x \mapsto x-y$ and let $\hat{Y}$ be the
composition of this map with the projection onto $P_{h(d)}(\mathcal{B}) / P_{h(d)+1}(\mathcal{B})$.  If $$F : (<\overline{b}>_{h(d)} \cap p^{-1} (\mathcal{B})_{h(d) +
2})/ <\overline{b}>_{h(d) + 1} \longrightarrow P_{h(d)}(\mathcal{B}) / P_{h(d)+1}(\mathcal{B})$$ is the map induced by $\hat{Y}$ on the quotient, then the following are 
equivalent:
\begin{enumerate} 
\item The range of $F$ is not all of $P_{h(d)}(\mathcal{B}) / P_{h(d)+1}(\mathcal{B})$. 
\item There exists in $P_{h(d)}(\mathcal{B})$ an element of height $h(d)$ which is proper with respect to $<\overline{b}>$.
\end{enumerate} 
\end{claim}

\begin{proof} To show 2 $\rightarrow$ 1, suppose $w \in P_{h(d)}$ has height $h(d)$ and is proper with respect to $<\overline{b}>$.  
Then the coset of $w$ is not in the range of $F$.  Otherwise, $w = x-y+q$ for some $x \in <\overline{b}>$, some $y \in (\mathcal{B})_{h(d)}$, and some $q \in
P_{h(d)+1}(\mathcal{B})$.  But then $h(w-x)>h(d)$, so $w$ was not proper.

To show the other implication, suppose that $w$ is an element of $P_{h(d)}(\mathcal{B})$ representing a coset not in the range of $F$.  Then $h(w) = h(d)$.  Further, $w$ 
is proper, since if it were not, and if $h(s-w)> h(d)$ witnessed this, we could write $s-w = p \zeta$ with $\zeta \in (\mathcal{B})_{h(d)}$.  But then $ps = p \zeta$ 
since $pw = 0$.  But then $F$ will map $s$ to the coset of $v$, giving a contradiction.
\end{proof}

Now since $d-v$ is such an element as is described in the second condition of the claim, we know that the range of $F$ is not all of $P_{h(d)}(\mathcal{B}) /
P_{h(d)+1}(\mathcal{B})$.  Since the vector spaces are finite (and thus finite dimensional), we know that the dimension of $(<\overline{b}>_{h(d)} \cap p^{-1}
(\mathcal{B})_{h(d) + 2})/ <\overline{b}>_{h(d) + 1}$ is less that $u_{h(d)}(\mathcal{B})$.  However, since $f$ was height preserving, it maps 
\[\begin{array}{c} (<\overline{b}>_{h(d)} \cap p^{-1} (\mathcal{B})_{h(d) + 2})/ <\overline{b}>_{h(d) + 1} \\
\downarrow \mbox{onto} \\ 
(<\overline{a}>_{h(d)} \cap p^{-1} (\mathcal{A})_{h(d) + 2})/ <\overline{a}>_{h(d) + 1}\\ \end{array}\]
Thus the dimension of $(<\overline{a}>_{h(d)} \cap p^{-1} (\mathcal{A})_{h(d) + 2})/ <\overline{a}>_{h(d) + 1}$ is less than
$u_{h(d)}(\mathcal{B})$.

In the case that $h(d)<\omega \cdot \delta + \omega$, we now know that the dimension of 
\[(~<~\overline{a}~>_{h(d)}~\cap~p^{-1}~(\mathcal{A})_{h(d) + 2})/
<\overline{a}>_{h(d) + 1}\] is less than $u_{h(d)}(\mathcal{A})$, so there 
is an element $c_1$ in $\mathcal{A}$ such that $pc_1 = 0$, $h(pc_1) = h(d)$, and which is 
proper with respect to $<\overline{a}>$.  Since $h(z) > h(d) + 1$, we may write $z = pc_2$ where $c_2 \in (\mathcal{B})_{h(d)+1}$.  Now we write $c = c_1 + c_2$ and note 
that $pc = z$, that $h(c) = h(d)$, and finally that $c$ is proper with respect to $<\overline{a}>$.

If $h(d) \geq \omega \cdot \delta + \omega$, we need considerably less.  In particular, it suffices to find some $c$ such that $pc=z$, such that $c$ is proper with respect
to $<\overline{a}>$, and such that $h(c) = \omega \cdot \delta + \omega$.  This can be achieved by replacing $h(d)$ with $\omega \cdot \delta + \omega$ in the preceding 
argument, and noting that since $\omega \cdot \delta$ is a limit, $u_{\omega \cdot \delta} = \omega$.  This completes the proof for the case $(\mathcal{A}, 
\overline{a}) \leq_{2 \cdot \delta} (\mathcal{B}, \overline{b})$ with $\delta$ a successor.

If $\delta$ is a limit ordinal, it suffices to consider some odd successor ordinal $2 \cdot \eta + 1 < 2 \cdot \delta$ and to show that for any $\overline{d} \in
\mathcal{B}$ there is some $\overline{c} \in \mathcal{A}$ such that $(\mathcal{B}, \overline{b}, \overline{d}) \leq_{2 \cdot \eta + 1} (\mathcal{A},
\overline{a}, \overline{c})$.  Then the proof is exactly as in the successor case.

In the case that we start with $(\mathcal{A}, \overline{a}) \leq_{2 \cdot \delta + 1} (\mathcal{B}, \overline{b})$ , we need to show that for any $\overline{d} \in
\mathcal{B}$ there is some $\overline{c} \in \mathcal{A}$ such that $(\mathcal{B}, \overline{b}, \overline{d}) \leq_{2 \cdot \delta} (\mathcal{A}, \overline{a},
\overline{c})$.  Now we can follow the proof exactly as in the even successor case, except that we replace $\omega \cdot \delta + \omega$ with $\omega \cdot \delta$.
\end{proof}

We now adapt the relations $\leq_\beta$ on pairs $(\mathcal{A}, \bar{a}), (\mathcal{B}, \bar{b})$ to relations on $L$.

\begin{dfn} We say that $(j_1, p_1) \leq_\beta (j_2, p_2)$ if and only if \[(G^{j_1}, ran(p_1)) \leq_\beta (G^{j_2}, ran(p_2))\]\end{dfn}

We need to verify that $(L, U, P, \hat{\ell}, E, (\leq_{\beta})_{\beta<\hat{\alpha}})$ is an $\hat{\alpha}$-system.  
For the necessary effectiveness, notice that we
need only consider $\leq_\beta$ on members of $L$, so only the groups $G^i$ are considered.  Conditions 1 -- 3 are clear, as is the fact that $(\leq_{\beta})_{\beta <
\hat{\alpha}}$ is uniformly c.e.  It remains to verify the following:

\begin{lem} If $\sigma u \in P$ where $\sigma$ ends in $\ell^0$ and $$\ell^0 \leq_{\beta_0} \ell^1 \leq_{\beta_1} \dots \leq_{\beta_{k-1}} \ell^k$$ where $\beta_0 > \beta_1
> \dots > \beta_k$, then there exists some $\ell^*$ such that $\sigma u \ell^* \in P$ and for all $i \leq k$, we have $\ell^i \leq_{\beta_0} \ell^*$.  \end{lem}

\begin{proof} We write $\ell^i = (j_i, p_i)$.  By Lemma \ref{baf}, given $\ell^{k-1} \leq_{\beta_{k-1}} \ell^k$ we can produce an $\tilde{\ell}^{k-1} = (\tilde{j}_{k-1}, \tilde{p}_{k-1})$ such that
$\tilde{p}$ extends $p_{k-1}$ (mapping into the same structure) and $\ell^k \leq_{\beta_k} \tilde{\ell}^{k-1}$.  Similarly, for each $i$, produce $\tilde{\ell}^i$ such that $\ell^{i+1} \leq_{\beta_{i+1}}
\tilde{\ell}^i$.  It will then be the case that for all $i$, $\ell^i \leq_{\beta_i} \tilde{\ell}^0$.  If $u = 0$ or if $1$ occurs somewhere in $\sigma$, let $\ell^* =
(\tilde{j}_0, p^*)$, where $p^*$ extends $\tilde{p}_0$ and its domain and range each contain the first $n$ constants, where $2n+1$ is the length of $\sigma$.  Now $\sigma u
\ell^* \in P$ and for all $i$, $\ell^i \leq_{\beta_0} \ell^*$.

If, on the other hand, $u = 1$ and $1$ does not occur in $\sigma$, then we may be sure that $\tilde{j}_0 = 0$.  In
this case, find some $j^*>0$ such that $\alpha_{j^*} > \beta_0$.  Note that since for each $\beta < \alpha_{j^*}$ we
have $u_\beta (G^{j^*}) = u_\beta (G^0)$, it follows that $(G^{j^*}, \emptyset) \leq_{\beta_0 + 1} (G^0, \emptyset)$.  
Thus, by Lemma \ref{baf}, we have some sequence $ran(p^*) \subseteq G^{j^*}$ such that $(G^0, ran(\tilde{p}))
\leq_{\beta_0} (G^{j^*}, ran(p^*))$ and having length $n$ where $2n+1$ is the length of $\sigma$.  We define $p^*$ to
be the function taking each of an initial sequence of the natural numbers to the corresponding element of that
sequence. Then clearly $\sigma u \ell^* \in P$, and for any $i$, we have $\ell^i \leq_{\beta_0} (G^0, ran(\tilde{p}))
\leq_{\beta_0} (G^{j^*}, ran(p^*))$ \end{proof}

Now let $S$ be an arbitrary $\Pi^0_{\hat{\alpha}}$ set.  There is a $\Delta^0_{\hat{\alpha}}$ function $g(n, s):  
\omega^2 \to 2$ such that for all $n$, we have $n \in S$ if and only if $\forall s [g(n, s) = 0]$, and such that for
all $n, s \in \omega$, if $g(n, s) = 1$ then $g(n, s+1) = 1$.  We define a $\Delta^0_{\hat{\alpha}}$ instruction
function $q_n$ as follows.  If $\sigma \in P$ and $\sigma$ is of length $m$, then we define $q_n(\sigma) = g(n, m)$.

Now we certainly can find computable indices for all the components of the $\hat{\alpha}$-system, and we can uniformly find a $\Delta^0_{\hat{\alpha}}$ index for each
$q_n$, so the Ash metatheorem gives us (uniformly in $n$), a run $\pi_n$ of $(P, q_n)$ and the index for the c.e.\ set 
$\bigcup\limits_{i \in \omega} E(\pi_n(2i))$.  Let
$H^n$ denote the group whose diagram this is.  Note that if $n \in S$, then $q_n (m) = 0$ for all $m$, and so $H^n \simeq G^0$.  Otherwise there is some $\hat{m}$ such that
for all $m > \hat{m}$, we have $q_n(m) = 1$, and so $H^n \simeq G^i$ for some $i \neq 0$.\end{proof}

\bibliographystyle{plain}
\bibliography{ecl}

\begin{thebibliography}{10}

\bibitem{ashlab}
C.~J. Ash.
\newblock Labelling systems and r.e.\ structures.
\newblock {\em Annals of Pure and Applied Logic}, 47:99--119, 1990.

\bibitem{akbook}
C.~J. Ash and J.~F. Knight.
\newblock {\em Computable structures and the hyperarithmetical hierarchy}.
\newblock Elsevier, 2000.

\bibitem{barker}
E.~Barker.
\newblock Back and forth relations for reduced {A}belian $p$-groups.
\newblock {\em Annals of Pure and Applied Logic}, 75:223--249, 1995.

\bibitem{ecl1}
W.~Calvert.
\newblock The isomorphism problem for classes of computable fields.
\newblock preprint, 2003.

\bibitem{frsimsm}
H.~Friedman, S.~Simpson, and R.~Smith.
\newblock Countable algebra and set existence axioms.
\newblock {\em Annals of Pure and Applied Logic}, 25:141--181, 1983.

\bibitem{fs}
H.~Friedman and L.~Stanley.
\newblock A {B}orel reducibility theory for classes of countable structures.
\newblock {\em Journal of Symbolic Logic}, 54:894--914, 1989.

\bibitem{gk}
S.~S. Goncharov and J.~F. Knight.
\newblock Computable structure and non-structure theorems.
\newblock {\em Algebra and Logic}, 41:351--373, 2002.

\bibitem{hkss}
D.~Hirschfeldt, B.~Khoussainov, R.~Shore, and A.~M. Slinko.
\newblock Degree spectra and computable dimensions in algebraic structures.
\newblock {\em Annals of Pure and Applied Logic}, 115:71--113, 2002.

\bibitem{hjbook}
G.~Hjorth.
\newblock {\em Classification and orbit equivalence relations}.
\newblock American Mathematical Society, 1999.

\bibitem{kaplansky}
I.~Kaplansky.
\newblock {\em Infinite {A}belian groups}.
\newblock University of Michigan Press, 1969.

\bibitem{khis}
N.~G. Khisamiev.
\newblock Constructive {A}belian $p$-groups.
\newblock {\em Siberian Advances in Mathematics}, 2:68--113, 1992.

\bibitem{clin}
C.~Lin.
\newblock The effective content of {U}lm's theorem.
\newblock In {\em Aspects of effective algebra}, pages 147--160. Upside Down A
  Book Company, 1979.

\bibitem{linjsl}
C.~Lin.
\newblock Recursively presented {A}belian groups: effective $p$-group theory
  {I}.
\newblock {\em Journal of Symbolic Logic}, 46:617 -- 624, 1981.

\bibitem{moro}
A.~S. Morozov.
\newblock Functional trees and automorphisms of models.
\newblock {\em Algebra and Logic}, 32:28--38, 1993.

\bibitem{nies}
A.~Nies.
\newblock Undecidable fragments of elementary theories.
\newblock {\em Algebra Universalis}, 35:8--33, 1996.

\bibitem{rabsco}
M.~O. Rabin and D.~Scott.
\newblock The undecidability of some simple theories.
\newblock preprint.

\bibitem{richman}
F.~Richman.
\newblock The constructive theory of countable {A}belian $p$-groups.
\newblock {\em Pacific Journal of Mathematics}, 45:621 -- 624, 1973.

\bibitem{lrogers}
L.~Rogers.
\newblock The structure of $p$-trees: algebraic systems related to {A}belian
  groups.
\newblock In {\em {A}belian Group Theory: 2nd {N}ew {M}exico State Conference},
  pages 57--72. Springer-Verlag, 1976.

\bibitem{simpson}
S.~Simpson.
\newblock {\em Subsystems of second order arithmetic}.
\newblock Springer-Verlag, 1999.

\end{thebibliography}

\end{document}